\documentclass[11pt]{amsart}
\newtheorem{theorem}{Theorem}[section]
\newtheorem{corollary}[theorem]{Corollary}
\newtheorem{lemma}[theorem]{Lemma}
\newtheorem{proposition}[theorem]{Proposition}

\theoremstyle{definition}

\newtheorem{remark}[theorem]{Remark}
\numberwithin{equation}{section}

\newcommand{\real}{\mathbb R}

\def\natu{\mathbb N}

\begin{document}

\title[Sharp estimates for
semi-stable radial solutions] {Sharp estimates for semi-stable
radial solutions of semilinear elliptic equations}
\author{Salvador Villegas}
\thanks{The author has been supported by the MEC Spanish grants
MTM2005-01331 and MTM2006-09282}
\address{Departamento de An\'{a}lisis
Matem\'{a}tico, Universidad de Granada, 18071 Granada, Spain.}
\email{svillega@ugr.es}

\begin{abstract}
This paper is devoted to the study of semi-stable radial solutions
$u\in H^1(B_1)$ of $-\Delta u=g(u) \mbox{ in } B_1\setminus \{
0\}$, where $g\in C^1(\real)$ is a general nonlinearity and $B_1$
is the unit ball of $\real^N$. We establish sharp pointwise
estimates for such solutions. As an application of these results,
we obtain optimal pointwise estimates for the extremal solution
and its derivatives (up to order three) of the semilinear elliptic
equation $-\Delta u=\lambda f(u)$, posed in $B_1$, with Dirichlet
data $u|_{\partial B_1}=0$, and a continuous, positive,
nondecreasing and convex function $f$ on $[0,\infty)$ such that
$f(s)/s\rightarrow\infty$ as $s\rightarrow\infty$.

In addition, we provide, for $N\geq 10$, a large family of
semi-stable radially decreasing unbounded $H^1(B_1)$ solutions.

\end{abstract}

\maketitle
\section{Introduction and main results}

This paper deals with the semi-stability of radial solutions $u\in
H^1(B_1)$ of
\begin{equation}\label{mainequation}
-\Delta u=g(u) \ \ \mbox{ in } B_1\setminus \{ 0\}\, ,
\end{equation}
\noindent where $B_1$ is the unit ball of $\real^N$, and $g\in
C^1(\real)$ is a general nonlinearity.

A radial solution $u\in H^1(B_1)$ of (\ref{mainequation}) is
called semi-stable if
$$\int_{B_1} \left( \vert \nabla v\vert^2-g'(u)v^2\right) \,
dx\geq 0$$ \noindent for every $v\in C^\infty (B_1)$ with compact
support in $B_1\setminus \{ 0\}$.

As an application of some general results obtained in this paper
for this class of solutions (for arbitrary $g\in C^1(\real)$), we
will establish sharp pointwise estimates related to the following
semilinear elliptic equation, which has been extensively studied.
$$
\left\{
\begin{array}{ll}
-\Delta u=\lambda f(u)\ \ \ \ \ \ \  & \mbox{ in } \Omega \, ,\\
u\geq 0 & \mbox{ in } \Omega \, ,\\
u=0  & \mbox{ on } \partial\Omega \, ,\\
\end{array}
\right. \eqno{(P_\lambda)}
$$
\

\noindent where $\Omega\subset\real^N$ is a smooth bounded domain,
$N\geq 2$, $\lambda\geq 0$ is a real parameter, and the
nonlinearity $f:[0,\infty)\rightarrow \real$ satisfies
\begin{equation}\label{convexa}
f \mbox{ is } C^1, \mbox{ nondecreasing and convex, }f(0)>0,\mbox{
and }\lim_{u\to +\infty}\frac{f(u)}{u}=+\infty.
\end{equation}

\

It is well known that there exists a finite positive extremal
parameter $\lambda^\ast$ such that ($P_\lambda$) has a minimal
classical solution $u_\lambda\in C^2(\overline{\Omega})$ if $0\leq
\lambda <\lambda^\ast$, while no solution exists, even in the weak
sense, for $\lambda>\lambda^\ast$. The set $\{u_\lambda:\, 0\leq
\lambda < \lambda^\ast\}$ forms a branch of classical solutions
increasing in $\lambda$. Its increasing pointwise limit
$u^\ast(x):=\lim_{\lambda\uparrow\lambda^\ast}u_\lambda(x)$ is a
weak solution of ($P_\lambda$) for $\lambda=\lambda^\ast$, which
is called the extremal solution of ($P_\lambda$) (see
\cite{Bre,BV}).

The regularity and properties of the extremal solutions depend
strongly on the dimension $N$, domain $\Omega$ and nonlinearity
$f$. When $f(u)=e^u$, it is known that $u^\ast\in L^\infty
(\Omega)$ if $N<10$ (for every $\Omega$) (see \cite{CrR,MP}),
while $u^\ast (x)=-2\log \vert x\vert$ and $\lambda^\ast=2(N-2)$
if $N\geq 10$ and $\Omega=B_1$ (see \cite{JL}). There is an
analogous result for $f(u)=(1+u)^p$ with $p>1$ (see \cite{BV}).
Brezis and V\'azquez \cite{BV} raised the question of determining
the boundedness of $u^\ast$, depending on the dimension $N$, for
general nonlinearities $f$ satisfying (\ref{convexa}). The best
result is due to Nedev \cite{Ne}, who proved that $u^\ast \in
L^\infty (\Omega)$ if $N\leq 3$, and Cabr\'e \cite{cabre4}, who
has proved recently that $u^\ast \in L^\infty (\Omega)$ if $N=4$
and $\Omega$ is convex. Cabr\'e and Capella \cite{cc} have proved
that $u^\ast \in L^\infty (\Omega)$ if $N\leq 9$ and $\Omega=B_1$
(similar results for the $p-$laplacian operator are contained in
\cite{plaplaciano}). Another interesting question is whether the
extremal solution lies in the energy class.  Nedev \cite{Ne,Ne2}
proved that $u^\ast \in H_0^1(\Omega)$ if $N\leq 5$ (for every
$\Omega$) or $\Omega$ is strictly convex (for every $N\geq 2$).
Brezis and V\'azquez \cite{BV} proved that a sufficient condition
to have $u^\ast \in H_0^1(\Omega)$ is that $\liminf_{u\to \infty}
u\, f'(u)/f(u)>1$ (for every $\Omega$ and $N\geq 2$). On the other
hand, it is an open problem (see \cite[Problem 5]{BV}) to know the
behavior of $f'(u^\ast)$ near the the singularities of $u^\ast$.
Is it always like $C/\vert x\vert^2\, $?

If $\Omega=B_1$, it is easily seen by the Gidas-Ni-Nirenberg
symmetry result that $u_\lambda$ is radially decreasing for
$0<\lambda<\lambda^\ast$. Hence, its limit $u^\ast$ is also
radially decreasing. In this situation, Cabr\'e and Capella
\cite{cc} have proved the following result:

\begin{theorem}(\cite{cc}).\label{cabrecapella}
Assume that $\Omega=B_1$, $N\geq 2$, and that $f$ satisfies
(\ref{convexa}). Let $u^\ast$ be the extremal solution of
($P_\lambda$). We have that

\begin{enumerate}

\item[i)] If $N<10$, then $u^\ast \in L^\infty (B_1)$,

\

\item[ii)] If $N=10$, then $u^\ast(x)\leq C\, \left\vert \log
\vert x\vert \right\vert $ \ in $B_1$ for some constant $C$,

\

\item[iii)] If $N>10$, then $\displaystyle{u^\ast(x)\leq C\, \vert
x\vert^{-N/2+\sqrt{N-1}+2}\sqrt{\left\vert \log \vert x\vert
\right\vert}} \ $ in $B_1$ for some constant $C$,

\

\item[iv)] If $N\geq 10$ and $k\in \{1,2,3\}$, then
$\displaystyle{\vert
\partial^{k}u^\ast(x)\vert \leq C\,
\,\vert x\vert^{-N/2+\sqrt{N-1}+2-k}\sqrt{\left\vert \log \vert
x\vert \right\vert}}\ $ in $B_1$ for some constant $C$.

\end{enumerate}

\end{theorem}

\

Among other results, in this paper we establish sharp pointwise
estimates for $u^\ast$ and its derivatives (up to order three) in
the radial case. We improve the above theorem, answering
affirmatively to an open question raised in \cite{cc}, about the
removal of the factor $\sqrt{\left\vert \log \vert x\vert
\right\vert}$.

By abuse of notation, we write $u(r)$ instead of $u(x)$, where
$r=\vert x\vert$ and $x\in \real^N$. We denote by $u_r$ the radial
derivative of a radial function $u$.

\begin{theorem}\label{extremal}

Assume that $\Omega=B_1$, $N\geq 2$, and that $f$ satisfies
(\ref{convexa}). Let $u^\ast$ be the extremal solution of
($P_\lambda$). We have that

\begin{enumerate}

\item[i)] If $N<10$, then $u^\ast(r)\leq C\, (1-r)\, , \ \ \forall
r\in [0,1]$,

\

\item[ii)] If $N=10$, then $u^\ast(r)\leq C\, \vert \log r\vert \,
, \ \ \forall r\in (0,1]$,

\

\item[iii)] If $N>10$, then $\displaystyle{u^\ast(r)\leq C\,
\left( r^{-N/2+\sqrt{N-1}+2}-1\right) \, , \ \ \forall r\in
(0,1]}$,

\

\item[iv)] If $N\geq 10$, then $\displaystyle{\vert
\partial_r^{(k)}u^\ast(r)\vert \leq C\,
\,r^{-N/2+\sqrt{N-1}+2-k} \, , \ \ \forall r\in (0,1]},$

\noindent $\forall k\in \{1,2,3\}$,

\end{enumerate}

\noindent where $\displaystyle{C=C_N \min_{t\in [1/2,1]}\vert
u^\ast_r(t)\vert}$, and $C_N$ is a constant depending only on $N$.

\end{theorem}

\,

\begin{remark}\label{C} It is immediate that if we replace the
function $f$ by $\tilde{f}:=f(\cdot/M)$, with $M>0$, then the
extremal solution $\tilde{u}^\ast$ associated to $\tilde{f}$ is
$\tilde{u}^\ast=M u^\ast$. Hence the constant $C$ in Theorem
\ref{extremal} must depend homogeneously on $u^\ast$. In fact,
this linear coefficient is very small since, for instance, we have
$$\min_{t\in [1/2,1]}\vert u^\ast_r(t)\vert\leq
4(u^\ast(1/2)-u^\ast(3/4))\leq 4u^\ast(1/2)\leq
\frac{4}{\mbox{measure}\, ( B_{1/2})} \Vert u^\ast
\Vert_{L^1(B_{1/2})}.$$
\end{remark}

\begin{remark} In \cite{BV} it is proved that if

$$N>10 \ \ \ \ \ \ \mbox{      and        }\ \ \ \ \ \ p\geq
p_N:=\frac{N-2\sqrt{N-1}}{N-2\sqrt{N-1}-4}\, ,$$

\noindent then the extremal solution for $f(u)=(1+u)^p$ and
$\Omega=B_1$ is given by $u^\ast(r)=r^{-2/(p-1)}-1$. In
particular, if $N>10$ and $p=p_N$ (called the Joseph-Lundgren
exponent), then $u^\ast(r)=r^{-N/2+\sqrt{N-1}+2}-1$. Hence the
pointwise estimates of Theorem \ref{extremal} for $u^\ast$ and its
derivatives (up to order three) are optimal if $N>10$. The
optimality of the theorem for $N=10$ follows immediately by
considering $f(u)=e^u$. As mentioned before, it is obtained in
this case that $u^\ast(r)=2\vert \log r\vert$.

\end{remark}

\begin{remark}

In fact, the convexity of $f$ is not necessary to obtain our main
results. Specifically, if we assume $f\in C^1$, nondecreasing,
$f(0)>0$ and $\lim_{u\to +\infty}f(u)/u=+\infty$, then it can be
proved (see \cite[Proposition 5.1]{cc}) that there exits a finite
positive extremal parameter $\lambda^\ast$ such that ($P_\lambda$)
has a minimal classical solution $u_\lambda\in
C^2(\overline{\Omega})$ if $0\leq \lambda <\lambda^\ast$, while no
solution exists, even in the weak sense, for
$\lambda>\lambda^\ast$. The set $\{u_\lambda:\, 0\leq \lambda <
\lambda^\ast\}$ of classical solutions is increasing in $\lambda$
and its  pointwise limit
$u^\ast(x):=\lim_{\lambda\uparrow\lambda^\ast}u_\lambda(x)$ is a
semi-stable weak solution of ($P_\lambda$) for
$\lambda=\lambda^\ast$. Note that the family of minimal solutions
$\{ u_\lambda \}$ may not be continuous as a function of
$\lambda$, as in the case of $f$ convex. Under these hypothesis of
$f$ it is possible to obtain the results (with the only exception
of the case $N\geq 10$ and $k=3$ of item iv)) of Theorems
\ref{cabrecapella} and \ref{extremal}.

\end{remark}

As we have mentioned, the proof of Theorem \ref{extremal} is based
on general properties of semi-stable radial solutions. Note that
the minimality of $u_\lambda$ implies  its semi-stability.
Clearly, we can pass to the limit and obtain that $u^\ast$ is also
radial and semi-stable. In addition, by a result of Nedev
\cite{Ne2} (see also \cite{cc}), we have that $u^\ast\in
H_0^1(B_1)$.

Recalling the definition of the semi-stability at the beginning of
the paper, we observe that a radial solution $u\in H^1(B_1)$ of
(\ref{mainequation}) is bounded away from the origin. Hence, using
standard regularity results, we obtain $u\in C^2(B_1\setminus \{
0\})$, and the definition of semi-stability makes sense.

If $u$ is a bounded radial solution of (\ref{mainequation}), then
$u\in C^2(\overline{B_1})$ and the semi-stability of $u$ means
that the first eigenvalue of the linearized problem
$-\Delta-g'(u)$ in $B_1$ is nonnegative.

Note that the expression which defines the semi-stability is
nothing but the second variation of the energy functional
associated to (\ref{mainequation}) in a domain
$\Omega\subset\real^N$ (with $\overline{\Omega}\subset
B_1\setminus\{ 0\}$): $E_\Omega (u)=\int_\Omega \left( \vert
\nabla u\vert^2 /2-G(u)\right) \, dx$, where $G'=g$. Thus, if
$u\in C^2(B_1\setminus \{ 0\})$ is a local minimizer of $E_\Omega$
for every smooth domain $\Omega\subset\real^N$ (with
$\overline{\Omega}\subset B_1\setminus\{ 0\}$) (i.e., a minimizer
under every small enough $C^1(\Omega)$ perturbation vanishing on
$\partial \Omega$), then $u$ is a semi-stable solution of
(\ref{mainequation}). Other general situations include stable
solutions: minimal solutions, extremal solutions or absolute
minimizers between a subsolution and a supersolution (see
\cite[Rem. 1.11]{cc} for more details).

Our main results about semi-stable radial solutions are the
following.

\begin{theorem}\label{principal}

Let $N\geq 2$, $g\in C^1(\real)$, and $u\in H^1(B_1)$ be a
semi-stable radial solution of (\ref{mainequation}). Then there
exists a constant $M_N$ depending only on $N$ such that:

\begin{enumerate}

\item[i)] If $N<10$, then $\Vert u\Vert_{L^\infty(B_1)}\leq M_N
\Vert u\Vert_{H^1(B_1\setminus \overline{B_{1/2}})}$.

\

\item[ii)] If $N=10$, then $\vert u(r)\vert \leq M_{10} \Vert
u\Vert_{H^1(B_1\setminus \overline{B_{1/2}})} \, (\vert \log
r\vert +1)\, , \ \ \forall r\in (0,1]$.

\

\item[iii)] If $N>10$, then $\displaystyle{\vert u(r)\vert \leq
M_N \Vert u\Vert_{H^1(B_1\setminus \overline{B_{1/2}})} \,
r^{-N/2+\sqrt{N-1}+2}\, , \ \ \forall r\in (0,1]}$.
\end{enumerate}

\end{theorem}

\

\begin{theorem}\label{estimas}

Let $N\geq 2$, $g\in C^1(\real)$, and $u\in H^1(B_1)$ be a
semi-stable radially decreasing solution of (\ref{mainequation}).
Then there exists a constant $M'_N$ depending only on $N$ such
that:

\begin{enumerate}

\item[i)] If $g\geq 0$, then

$$\vert u_r(r)\vert \leq M'_N \Vert \nabla
u\Vert_{L^2(B_1\setminus B_{1/2})} r^{-N/2+\sqrt{N-1}+1}\, ,\ \
\forall r\in (0,1/2].$$

\

\item[ii)] If $g\geq 0$ is nondecreasing, then

$$\vert u_{rr}(r)\vert \leq M'_N \Vert \nabla
u\Vert_{L^2(B_1\setminus B_{1/2})} r^{-N/2+\sqrt{N-1}}\, ,\ \
\forall r\in (0,1/2].$$

\

\item[iii)] If $g\geq 0$ is nondecreasing and convex, then

$$\vert
u_{rrr}(r)\vert \leq M'_N \Vert \nabla u\Vert_{L^2(B_1\setminus
B_{1/2})} r^{-N/2+\sqrt{N-1}-1}\, ,\ \ \forall r\in (0,1/2].$$

\end{enumerate}

\end{theorem}

\

\begin{remark}

\label{anillo} We emphasize that the estimates obtained in
Theorems \ref{principal} and \ref{estimas} are in terms of the
$H^1$ norm of the annulus $B_1\setminus \overline{B_{1/2}}$, while
$u$ is required to belong to $H^1(B_1)$. In fact, this requirement
is essential to obtain our results, since we can always find
radial weak solutions of (\ref{mainequation}) (not in the Sobolev
space of the unit ball), for which the statements of Theorems
\ref{principal} and \ref{estimas} fail to satisfy (see
\cite{BV,cc}).

\end{remark}

\begin{remark}

In \cite[Rem. 1.9]{cc} the authors raised the question whether the
estimates of Theorem \ref{estimas} hold for general nonlinearities
$g$, without the assumptions on the nonnegativeness of $g$, $g'$
and/or $g''$. In this paper we answer negatively to this question.
In fact, without assumptions on the sign of $g$, $g'$ or $g''$ it
is not possible to obtain any pointwise estimate for $\vert
u_r\vert$, $\vert u_{rr}\vert$ or $\vert u_{rrr}\vert$ (see
Corollaries \ref{nohay}, \ref{nohayy} and \ref{nohayyy}).

\end{remark}

To prove the main results of the paper we will use Lemma
\ref{essential}, which, roughly speaking, says that the are some
restrictions on the growth of the derivative of a radial
semi-stable solution of (\ref{mainequation}) around the origin. In
the proof of this lemma, we will make use of \cite[Lem. 2.1]{cc},
which was inspired by the proof of Simons theorem on the
nonexistence of singular minimal cones in $\real^N$ for $N\leq 7$
(see \cite[Th. 10.10]{simons} and \cite[Rem. 2.2]{cc} for more
details). Similar methods are used in \cite{cabre,yo} to study the
stability or instability of radial solutions in all space
$\real^N$.

\

The paper is organized as follows. In Section \ref{dos} we prove
Theorems \ref{extremal}, \ref{principal} and \ref{estimas}.
Section \ref{tres} provides, for $N\geq 10$, a large family of
semi-stable radially decreasing unbounded $H^1(B_1)$ solutions of
problems of the type (\ref{mainequation}). Taking solutions of
this family, we will show the impossibility of obtaining pointwise
estimates for $\vert u_r\vert$, $\vert u_{rr}\vert$ or $\vert
u_{rrr}\vert$ if no further assumptions on the sign of $g$, $g'$
or $g''$ are imposed.

\section{Proof of the main results}\label{dos}

\begin{lemma}\label{essential}

Let $N\geq 2$, $g\in C^1(\real)$, and $u\in H^1(B_1)$ be a
semi-stable radial solution of (\ref{mainequation}). Then there
exists a constant $K_N$ depending only on $N$ such that:

\begin{equation}\label{inequality}
\int_0^r t^{N-1}u_r(t)^2 \, dt\leq K_N \Vert \nabla
u\Vert_{L^2(B_1\setminus B_{1/2})}^2\, r^{2\sqrt{N-1}+2}\,
 \ \ \ \forall r\in [0,1].
\end{equation}

\end{lemma}

\noindent {\bf Proof.} Let us use \cite[Lem. 2.1]{cc} (see also
the proof of \cite[Lem. 2.3]{cc}) to assure that

$$(N-1)\int_{B_1}u_r^2\, \eta^2\, dx \leq
\int_{B_1}u_r^2\vert \nabla \left(r\, \eta\right) \vert^2 \, dx\,
,
$$

\noindent for every $\eta \in (H^1\cap L^\infty )(B_1)$ with
compact support in $B_1$ and such that $\vert \nabla \left(r\,
\eta\right) \vert\in L^\infty (B_1)$. Applying this inequality to
a radial function $\eta (\vert x\vert)$ we obtain

\begin{equation}\label{stablepropert}
(N-1)\int_0^1 u_r(t)^2\eta(t)^2 t^{N-1}\, dt\leq \int_0^1 u_r(t)^2
\left(t \, \eta (t)\right)'^{\, 2} t^{N-1}\, dt\, .
\end{equation}

We now fix $r\in (0,1/2)$ and consider the function

$$\eta (t)=\left\{
\begin{array}{ll}
r^{-\sqrt{N-1}-1} & \mbox{ if   } 0\leq t \leq r\, , \\ \\
t^{-\sqrt{N-1}-1} & \mbox{ if   } r<t\leq 1/2\, , \\ \\
2^{\sqrt{N-1}+2}(1-t) & \mbox{ if   } 1/2<t\leq 1\, .
\end{array}
\right.
$$

Since $(N-1)\eta(t)^2=\left(t \, \eta (t)\right)'^{\, 2}$ for
$r<t<1/2$, inequality (\ref{stablepropert}) shows that

$$\begin{array}{l}\displaystyle{\ \ \ (N-2)r^{-2\sqrt{N-1}-2}\int_0^r
u_r(t)^2 t^{N-1}\, dt}\\
\displaystyle{=\int_0^r \left( (N-1) \eta (t)^2-\left(t
\, \eta (t)\right)'^{\, 2}\right) u_r(t)^2 t^{N-1}\, dt }\\
\displaystyle{\leq -\int_{1/2}^1 \left( (N-1) \eta (t)^2-\left(t
\, \eta (t)\right)'^{\, 2}\right) u_r(t)^2 t^{N-1}\, dt }\leq
\displaystyle{ \alpha_N \int_{1/2}^1 u_r(t)^2 t^{N-1}\, dt,}
\end{array}$$

\noindent where the constant $\displaystyle{\alpha_N=\max_{1/2\leq
t\leq 1}-\left( (N-1) \eta (t)^2-\left(t \, \eta (t)\right)'^{\,
2}\right)}$ depends only on $N$. This establishes
(\ref{inequality}) for $r\in [0,1/2]$, if $N>2$.

If $r\in (1/2,1]$ and $N>2$ then, applying the above inequality
for $r=1/2$, we obtain

$$\begin{array}{l}\displaystyle{\ \ \ \int_0^r t^{N-1}u_r(t)^2\, dt\leq
\int_0^{1/2} t^{N-1}u_r(t)^2\, dt+\int_{1/2}^1 t^{N-1}u_r(t)^2\,
dt} \\
\displaystyle{\leq\left(
\frac{\alpha_N}{N-2}\left(\frac{1}{2}\right)^{2\sqrt{N-1}+2}+1\right)\int_{1/2}^1
t^{N-1}u_r(t)^2\, dt} \\
\displaystyle{\leq (2r)^{2\sqrt{N-1}+2}\left(
\frac{\alpha_N}{N-2}\left(\frac{1}{2}\right)^{2\sqrt{N-1}+2}+1\right)\int_{1/2}^1
t^{N-1}u_r(t)^2\, dt}
\end{array}$$

\noindent which is the desired conclusion with
$\displaystyle{K_N=\frac{1}{\omega_N}\left(\frac{\alpha_N}{N-2}+2^{2\sqrt{N-1}+2}\right)}$
(Note that the constant obtained for $r\in (1/2,1]$ is greater
than the one for $r\in [0,1/2]$).

Finally, if $N=2$, changing the definition of $\eta (t)$ in
$[0,r]$ by $\eta (t)=1/(r\, t)$, if $r_0<t\leq r$; $\eta(t)=1/(r\,
r_0)$, if $0\leq t\leq r_0$ (for arbitrary $r_0\in (0,r)$), we
obtain

$$\frac{1}{r^2}\int_{r_0}^r\frac{u_r(t)^2}{t}dt\leq \alpha_2
\int_{1/2}^1u_r(t)^2t\, dt\, .$$

Letting $r_0\rightarrow 0$ and taking into account that $t/r^2\leq
1/t$ for $0<t\leq r$ yields (\ref{inequality}) for $N=2$ and
$r\in[0,1/2]$. If $r\in (1/2,1]$, we can apply similar arguments
to the case $N>2$ to complete the proof. \qed

\begin{proposition}\label{prorand2r}

Let $N\geq 2$, $g\in C^1(\real)$, and $u\in H^1(B_1)$ be a
semi-stable radial solution of (\ref{mainequation}). Then there
exists a constant $K'_N$ depending only on $N$ such that:

\begin{equation}\label{inequalityrand2r}
\left\vert u(r)-u\left(\frac{r}{2}\right) \right\vert \leq K'_N
\Vert \nabla u\Vert_{L^2(B_1\setminus B_{1/2})}\,
r^{-N/2+\sqrt{N-1}+2}\,
 \ \ \ \forall r\in (0,1].
\end{equation}

\end{proposition}

\noindent {\bf Proof.} Fix $r\in (0,1]$. Applying Cauchy-Schwarz
and Lemma \ref{essential} we deduce

$$\begin{array}{l}
\displaystyle{\ \ \ \left\vert u(r)-u\left(\frac{r}{2}\right)
\right\vert \leq \int_{r/2}^r\vert
u_r(t)\vert t^\frac{N-1}{2}\frac{1}{t^\frac{N-1}{2}}\, dt} \\
\displaystyle{\leq\left( \int_{r/2}^r u_r(t)^2 t^{N-1}\,
dt\right)^{1/2} \left( \int_{r/2}^r \frac{1}{t^{N-1}}\,
dt\right)^{1/2}}\\
\displaystyle{\leq K_N^{1/2} \Vert \nabla u\Vert_{L^2(B_1\setminus
B_{1/2})}\, r^{\sqrt{N-1}+1}\left( r^{2-N} \int_{1/2}^1
\frac{1}{t^{N-1}}\, dt\right)^{1/2}},
\end{array}$$

\noindent and (\ref{inequalityrand2r}) is proved. \qed

\

\noindent {\bf Proof of Theorem \ref{principal}.} Let $0<r\leq 1$.
Then, there exist $m\in \natu$ and $1/2<r_1\leq 1$ such that
$r=r_1/2^{m-1}$. Since $u$ is radial we have $u(r_1)\leq \Vert
u\Vert_{L^\infty (B_1\setminus B_{1/2})}\leq \gamma_N \Vert
u\Vert_{H^1(B_1\setminus \overline{B_{1/2}})}$, where $\gamma_N$
depends only on $N$. From this and Proposition \ref{prorand2r}, it
follows that

\begin{equation}\label{key}
\begin{array}{l} \displaystyle{\vert u(r)\vert \leq \vert
u(r_1)-u(r)\vert +\vert u(r_1)\vert\leq \sum_{i=1}^{m-1}\left\vert
u\left(\frac{r_1}{2^{i-1}}\right)-u\left(\frac{r_1}{2^i}\right)\right\vert+\vert
u(r_1)\vert} \\
\displaystyle{\leq K'_N \Vert \nabla u\Vert_{L^2(B_1\setminus
B_{1/2})}
\sum_{i=1}^{m-1}\left(\frac{r_1}{2^{i-1}}\right)^{-N/2+\sqrt{N-1}+2}+
\gamma_N \Vert u\Vert_{H^1(B_1\setminus \overline{B_{1/2}})} }\\
\displaystyle{\leq \left(
K'_N\sum_{i=1}^{m-1}\left(\frac{r_1}{2^{i-1}}\right)^{-N/2+\sqrt{N-1}+2}+
\gamma_N \right) \Vert u\Vert_{H^1(B_1\setminus
\overline{B_{1/2}})}.}
\end{array}
\end{equation}

\

$\bullet$ If $2\leq N<10$, we have $-N/2+\sqrt{N-1}+2>0$. Then

$$\sum_{i=1}^{m-1}\left(\frac{r_1}{2^{i-1}}\right)^{-N/2+\sqrt{N-1}+2}\leq
\sum_{i=1}^{\infty}\left(\frac{1}{2^{i-1}}\right)^{-N/2+\sqrt{N-1}+2},$$

\noindent which is a convergent series. Applying (\ref{key}),
statement i) of the theorem is proved.

\

$\bullet$ If $N=10$, we have $-N/2+\sqrt{N-1}+2=0$. From
(\ref{key}) we obtain

$$\begin{array}{l} \displaystyle{\vert u(r)\vert \leq \left( K'_N
(m-1)+\gamma_N \right) \Vert u\Vert_{H^1(B_1\setminus
\overline{B_{1/2}})}}
\\=\displaystyle{\left( K'_N \left( \frac{\log r_1-\log
r}{\log 2}\right) + \gamma_N \right) \Vert
u\Vert_{H^1(B_1\setminus \overline{B_{1/2}})}}
\\ \displaystyle{\leq \left( \frac{K'_N}{\log 2}+ \gamma_N \right)
\left( \vert \log r\vert +1\right) \Vert u\Vert_{H^1(B_1\setminus
\overline{B_{1/2}})},}
\end{array}$$

\noindent which gives statement ii).

\

$\bullet$ If $N>10$, we have $-N/2+\sqrt{N-1}+2<0$. Then

$$\sum_{i=1}^{m-1}\left(\frac{r_1}{2^{i-1}}\right)^{-N/2+\sqrt{N-1}+2}=
\frac{r^{-N/2+\sqrt{N-1}+2}-r_1^
{-N/2+\sqrt{N-1}+2}}{(1/2)^{-N/2+\sqrt{N-1}+2}-1}.$$

From this and (\ref{key}), we conclude

$$\vert u(r)\vert \leq \left(
\frac{K'_N}{(1/2)^{-N/2+\sqrt{N-1}+2}-1}+\gamma_N \right)
r^{-N/2+\sqrt{N-1}+2}\Vert u\Vert_{H^1(B_1\setminus
\overline{B_{1/2}})}\, ,$$

\noindent which completes the proof. \qed

\

\noindent {\bf Proof of Theorem \ref{estimas}.}

\begin{enumerate}

\item[i)] We first observe that $(-r^{N-1}u_r)'=r^{N-1}g(u)\geq
0$. Hence $-r^{N-1}u_r$ is a positive nondecreasing function and
so is $r^{2N-2}u_r^2$. Thus, for $0<r\leq 1/2$, we have

$$\int_0^{2r} t^{N-1}u_r(t)^2\, dt\geq \int_r^{2r} t^{N-1}u_r(t)^2
\, dt=\int_r^{2r} t^{2N-2}u_r(t)^2 \frac{1}{t^{N-1}}\, dt$$

$$\geq r^{2N-2}u_r(r)^2 \int_r^{2r}\frac{1}{t^{N-1}}\,
dt=r^{2N-2}u_r(r)^2 \, r^{2-N}\int_1^{2}\frac{1}{t^{N-1}}\, dt\
,$$

From this and Lemma \ref{essential} we obtain i).

\

\item[ii)] Consider the function $\Psi(r)=-N\,
r^{1-1/N}u_r(r^{1/N})\, , r\in (0,1]$. It is easy to check that
$\Psi'(r)=g(u(r^{1/N}))\, , r\in (0,1]$. As $g$ is nonnegative and
nondecreasing we have that $\Psi$ is a nonnegative nondecreasing
concave function. It follows immediately that $0\leq \Psi'(r)\leq
\Psi(r)/r\, ,r\in (0,1]$; which becomes

$$0\leq
-(N-1)r^{-1/N}u_r(r^{1/N})-u_{rr}(r^{1/N})\leq -N\,
r^{-1/N}u_r(r^{1/N})\, ,r\in (0,1].$$

Hence
$$r^{-1/N}u_r(r^{1/N})\leq u_{rr}(r^{1/N})\leq
-(N-1)r^{-1/N}u_r(r^{1/N})\, ,r\in (0,1].$$

Therefore $\vert u_{rr}(r)\vert\leq (N-1)\vert u_r(r)\vert /r \,
,r\in (0,1]$; and ii) follows from i).

\

\item[iii)] An easy computation shows that

$$u_{rrr}=-u_r g'(u)-\frac{N-1}{r}u_{rr}+\frac{N-1}{r^2}u_r \, ,
\ \ r\in (0,1].$$

On the other hand, it is proved in \cite[Th. 1.8 (c)]{cc} that
$g'(u(r))\leq h_N/r^2\, , r\in (0,1]$, for some constant $h_N$.
Since we have shown $\vert u_{rr}(r)\vert\leq (N-1)\vert
u_r(r)\vert /r \, ,r\in (0,1]$ in the proof of statement ii), it
follows from the above formula $\vert u_{rrr}(r)\vert\leq s_N\vert
u_r(r)\vert /r^2 \, ,r\in (0,1]$, for some constant $s_N$
depending only on $N$. Recalling i), the proof is now completed.
\qed

\end{enumerate}

\

To deduce Theorem \ref{extremal} from Theorems \ref{principal} and
\ref{estimas} we need the following lemma.

\begin{lemma}\label{monotonias}

Let $N\geq 2$, $g\in C^1(\real)$ nonnegative and nondecreasing
function and $u$ a radially decreasing solution of
(\ref{mainequation}) (neither $u\in H^1(B_1)$ nor $u$ is
semi-stable is required). Then

\begin{enumerate}

\item[i)] $r^{N-1}\vert u_r\vert$ is nondecreasing for $r\in
(0,1]$.

\medskip

\item[ii)] $r^{-1}\vert u_r\vert$ is nonincreasing for $r\in
(0,1]$.

\medskip

\item[iii)] $\max_{t\in [1/2,1]}\vert u_r(t)\vert\leq
2^{N-1}\min_{t\in [1/2,1]}\vert u_r(t)\vert$.

\medskip

\item[iv)] $\Vert \nabla u\Vert_{L^2(B_1\setminus B_{1/2})}\leq
q_N \min_{t\in [1/2,1]}\vert u_r(t)\vert$, for a certain constant
$q_N$ depending only on $N$.

\end{enumerate}

\end{lemma}

\noindent {\bf Proof.} \begin{enumerate} \item[i)] Since $u_r<0$
we have $\left(r^{N-1}\vert u_r\vert\right)'=r^{N-1}g(u)\geq 0$.

\item[ii)] As in the proof of statement ii) of Theorem
\ref{estimas} we have that the function $\Psi(r)=-N\,
r^{1-1/N}u_r(r^{1/N})$ is nonnegative, nondecreasing and concave
for $r\in (0,1]$. Therefore $\Psi(r)/r=-N\, r^{-1/N}u_r(r^{1/N})$
is nonincreasing , and ii) follows immediately.

\item[iii)]Take $r_1,r_2 \in [1/2,1]$ such that $\vert
u_r(r_1)\vert=\min_{t\in [1/2,1]}\vert u_r(t)\vert$ and $\vert
u_r(r_2)\vert=\max_{t\in [1/2,1]}\vert u_r(t)\vert$.

If $r_2\leq r_1$, we deduce from i) that $\vert u_r(r_2)\vert \leq
(r_1/r_2)^{N-1}\vert u_r(r_1)\vert\leq 2^{N-1}\vert
u_r(r_1)\vert$.

If $r_2> r_1$, we deduce from ii) that $\vert u_r(r_2)\vert \leq
(r_2/r_1)\vert u_r(r_1)\vert\leq 2\vert u_r(r_1)\vert\leq
2^{N-1}\vert u_r(r_1)\vert$.

\medskip

\item[iv)] We see at once that
$$\Vert \nabla u\Vert_{L^2(B_1\setminus B_{1/2})}\leq (\mbox{measure}\,
(B_1\setminus B_{1/2}))^{1/2}\max_{t\in [1/2,1]}\vert
u_r(t)\vert\, ,$$ \noindent and iv) follows from iii). \qed

\end{enumerate}

\

\noindent {\bf Proof of Theorem \ref{extremal}}

As we have mentioned, it is well known that $u^\ast$ is a
semi-stable radially decreasing $H_0^1(B_1)$ solution of
(\ref{mainequation}) for $g(s)=\lambda^\ast f(s)$. Hence, we can
apply to $u^\ast$ the results obtained in Theorems \ref{principal}
and \ref{estimas} and Lemma \ref{monotonias}.

Let us first prove i), ii) and iii) for $r\in (0,1/2)$. Since
$u^\ast (1)=0$, and on account of statement iv) of Lemma
\ref{monotonias}, we have $\Vert u^\ast \Vert_{H^1(B_1\setminus
\overline{B_{1/2}})}\leq h_N \Vert \nabla u^\ast
\Vert_{L^2(B_1\setminus B_{1/2})}\leq h'_N \min_{t\in
[1/2,1]}\vert u^\ast_r(t)\vert$, for certain constants $h_N,h'_N$
depending only on $N$. From this and Theorem \ref{principal}:

\medskip

i) follows from the inequality $1\leq 2(1-r)$, for $r\in (0,1/2)$.

\medskip

ii) follows from the inequality $\displaystyle{\vert \log r\vert
+1\leq \frac{\log 2+1}{\log 2} \vert \log r\vert}$, for $r\in
(0,1/2)$.

\medskip

iii) follows from the inequality

$$r^{-N/2+\sqrt{N-1}+2}\leq
\frac{(1/2)^{-N/2+\sqrt{N-1}+2}}{(1/2)^{-N/2+\sqrt{N-1}+2}-1}(r^{-N/2+\sqrt{N-1}+2}-1),
\mbox{ for }r\in (0,1/2).$$

\

We next show i), ii) and iii) for $r\in [1/2,1]$. From statement
iii) of Lemma \ref{monotonias} it follows that

$$u^\ast (r)=\int_{r}^1 \vert u^\ast_r(t)\vert\, dt \leq (1-r)\,
2^{N-1}\min_{t\in [1/2,1]}\vert u^\ast_r(t)\vert \, , \ \ \ \
\forall r\in[1/2,1],$$

\noindent which is the desired conclusion if $N<10$. If $N=10$,
our claim follows from the inequality $1-r\leq \vert\log r\vert$,
for $r\in[1/2,1]$. Finally, if $N>10$, the desired conclusion
follows immediately from the inequality $1-r\leq z_N(
r^{-N/2+\sqrt{N-1}+2}-1)$, for $r\in [1/2,1]$, for a certain
constant $z_N$.

\

We now prove statement iv). In the case $k=1$ and $r\in (0,1/2)$,
it follows immediately from statement i) of Theorem \ref{estimas}
and statement iv) of Lemma \ref{monotonias}. The case $k=1$ and
$r\in [1/2,1]$ is also obvious on account of statement iii) of
Lemma \ref{monotonias} and the inequality $1\leq
r^{-N/2+\sqrt{N-1}+1}$, for $r\in [1/2,1]$, for $N\geq 10$.

Finally, as in the proof of statement ii) and iii) of Theorem
\ref{estimas}, we have $\vert u^\ast_{rr}(r)\vert\leq (N-1)\vert
u^\ast_r(r)\vert /r$ and $\vert u^\ast_{rrr}(r)\vert\leq s_N \vert
u^\ast_r(r)\vert /r^2$, for $r\in (0,1]$, which gives statement
iv) for $k=2,3$ from the case $k=1$. \qed

\section{A family of semi-stable solutions}\label{tres}

\begin{theorem}\label{familia}

Let $h\in (C^2\cap L^1)(0,1]$ be a nonnegative function and
consider

$$\Phi (r)=r^{2\sqrt{N-1}}\left( 1+\int_0^r h(s)\,
ds\right)\ \ \ \ \forall r\in (0,1].$$

Define $u_r<0$ by

$$\Phi'(r)=(N-1)\, r^{N-3}u_r(r)^2\ \ \ \ \forall r\in (0,1].$$

Then, for $N\geq 10$, $u$ is a semi-stable radially decreasing
unbounded $H^1(B_1)$ solution of a problem of the type
(\ref{mainequation}), where $u$ is any function with radial
derivative $u_r$.

\end{theorem}

To prove Theorem \ref{familia} we need the following lemma, which
is a generalization of the classical Hardy inequality:

\begin{lemma}\label{masquehardy}

Let $\Phi \in C^1(0,L)$, $0<L\leq \infty$, satisfying $\Phi'>0$.
Then

$$\int_0^L \frac{4 \Phi^2}{\Phi'} \xi'^2 \geq \int_0^L \Phi'
\xi^2\, ,$$

\noindent for every $\xi \in C^\infty (0,L)$ with compact support.

\end{lemma}

\noindent {\bf Proof.} Integrating by parts and applying
Cauchy-Schwarz we obtain

$$\int_0^L \Phi' \xi^2=-2 \int_0^L \Phi \xi \xi' \leq 2\int_0^L
\frac{\vert \Phi\vert}{\sqrt{\Phi'}} \vert \xi' \vert \sqrt{\Phi'}
\vert \xi \vert \leq 2\left(\int_0^L \frac{
\Phi^2}{\Phi'}\xi'^2\right)^{1/2}\left(\int_0^L \Phi'
\xi^2\right)^{1/2},$$

\noindent which establishes the desired inequality. \qed

\

In the case $\Phi(r)=((N-2)/4)r^{N-2}$, $r>0$, the above lemma is
the Hardy inequality for radial functions in $\real^N$, $N>2$.

\

\noindent {\bf Proof of Theorem \ref{familia}.} First of all,
since $\Phi \in C^1(0,1]\cap C[0,1]$ is an increasing function, we
obtain $\Phi'\in L^1(0,1)$ and hence $r^{N-1}u_r^2=r^2 \Phi'
/(N-1) \in L^1(0,1)$, which gives $u\in H^1(B_1)$.

On the other hand, since $\Phi'(r)\geq 2\sqrt{N-1}\,
r^{2\sqrt{N-1}-1},\, r\in(0,1]$, we deduce $\vert u_r(r)\vert \geq
\sqrt{2}(N-1)^{-1/4} \ r^{-N/2+\sqrt{N-1}+1},\, r\in (0,1]$. As
$N\geq 10$, we have $-N/2+\sqrt{N-1}+1\leq -1$. It follows that
$u_r\notin L^1(0,1)$ and, since $u$ is radially decreasing, we
obtain $\lim_{r\to 0}u(r)=+\infty$.

Since $h\in C^2(0,1]$,  it follows that $u_r \in C^2(0,1]$.
Therefore, $\Delta u\in C^1\left(\overline{B_1}\setminus\{
0\}\right)$. Hence, taking $g\in C^1(\real)$ such that
$g(s)=-\Delta u(u^{-1}(s))$, for $s\in[u(1),+\infty)$, we conclude
that $u$ is solution of a problem of the type
(\ref{mainequation}).

It remains to prove that $u$ is semi-stable. Taking into account
that $u_r\neq 0$ in $(0,1]$ and applying \cite[Lem. 2.1]{cc}, the
semi-stability of $u$ is equivalent to

\begin{equation}\label{stableproperty}
\int_0^1 r^{N-1}u_r^2\ \xi'^2\, dr \geq (N-1)\int_0^1 r^{N-3}
u_r^2\ \xi^2\, dr,
\end{equation}

\noindent for every $\xi \in C^\infty (0,1)$ with compact support.

For this purpose, we will apply the lemma above. From the
definition of $\Phi$ it is easily seen that $\Phi'\geq
2\sqrt{N-1}\, \Phi /r,$ $r\in (0,1]$. It follows that

$$\frac{\Phi'r^2}{N-1}\geq \frac{4 \Phi^2}{\Phi'}\mbox{  in
}(0,1].$$

Finally, since $\Phi'r^2/(N-1)=r^{N-1}u_r^2$ and
$\Phi'=(N-1)r^{N-3} u_r^2$ in $(0,1]$, we deduce
(\ref{stableproperty}) by applying Lemma \ref{masquehardy}. \qed

\

As an application of Theorem \ref{familia} we have the following
results, which show the impossibility of obtaining any pointwise
estimate for $\vert u_r\vert$, $\vert u_{rr}\vert$ or $\vert
u_{rrr}\vert$ if the positivity of $g$, $g'$ or $g''$ is not
satisfied, for semi-stable radially decreasing $H^1(B_1)$
solutions of a problem of the type (\ref{mainequation}) and $N\geq
10$.

\begin{proposition}\label{sucesiones}

Let $\{r_n\}\subset(0,1]$, $\{M_n\}\subset \real^+$ two sequences
with $r_n\downarrow 0$. Then, for $N\geq 10$, there exists $u\in
H^1(B_1)$, which is a semi-stable radially decreasing unbounded
solution of a problem of the type (\ref{mainequation}), satisfying

$$\vert u_r(r_n)\vert \geq M_n \ \ \ \ \forall n\in \mathbb{N}.$$

\end{proposition}

\noindent {\bf Proof.} It is easily seen that for every sequences
$\{r_n\}\subset(0,1]$, $\{y_n\}\subset \real^+$, with
$r_n\downarrow 0$, there exists a nonnegative function $h\in
(C^2\cap L^1)(0,1]$ satisfying $h(r_n)=y_n$. Take $y_n=(N-1)\,
M_n^2\, r_n^{N-2\sqrt{N-1}-3}$ and apply Theorem \ref{familia}
with this function $h$. It is clear, from the definition of
$\Phi$, that $\Phi'(r)\geq h(r) r^{2\sqrt{N-1}}, \, r\in (0,1]$.
Hence

$$(N-1)\, r_n^{N-3}u_r(r_n)^2=\Phi'(r_n)\geq h(r_n)
r_n^{2\sqrt{N-1}}=y_n r_n^{2\sqrt{N-1}}=(N-1)r_n^{N-3}M_n^2,$$

\noindent and the proposition follows. \qed

\begin{corollary}\label{nohay}

Let $N\geq 10$. There does not exist a function $\psi:(0,1]
\rightarrow \real^+$ with the following property: for every $u\in
H^1(B_1)$ semi-stable radially decreasing solution of a problem of
the type (\ref{mainequation}), there exist $C>0$ and $\varepsilon
\in (0,1]$ such that $\vert u_r(r)\vert \leq C \psi (r)$ for $r\in
(0,\varepsilon]$.

\end{corollary}

\noindent {\bf Proof.} Suppose that such a function $\psi$ exists
and consider the sequences $r_n=1/n$, $M_n =n \, \psi (1/n)$. By
the proposition above, there exists $u\in H^1(B_1)$, which is a
semi-stable radially decreasing unbounded solution of a problem of
the type (\ref{mainequation}), satisfying $\vert u_r(1/n)\vert
\geq n \, \psi (1/n)$, a contradiction. \qed

\begin{proposition}\label{sucesiones2}

Let $\{r_n\}\subset(0,1]$, $\{M_n\}\subset \real^+$ two sequences
with $r_n\downarrow 0$. Then, for $N\geq 10$, there exists $u\in
H^1(B_1)$, which is a semi-stable radially decreasing unbounded
solution of a problem of the type (\ref{mainequation}) with $g\geq
0$, satisfying

$$\vert u_{rr}(r_n)\vert \geq M_n \ \ \ \ \forall n\in \mathbb{N}.$$

\end{proposition}

\noindent {\bf Proof.} Let $h\in C^2(0,1]$, increasing, satisfying
$0\leq h\leq 1$. Define $\Phi$ and $u_r$ as in Theorem
\ref{familia}. We claim that

\begin{enumerate}
\item[i)] $u$ is a semi-stable radially decreasing unbounded
$H^1(B_1)$ solution of a problem of the type (\ref{mainequation})
with $g\geq 0$.

\item[ii)] $\vert u_r \vert\leq D_N r^{-N/2+\sqrt{N-1}+1} ,\ \
\forall r\in (0,1]$, where $D_N$ only depends on $N$.

\item[iii)] $-u_{rr}\geq E_N h'(r)r^{-N/2+\sqrt{N-1}+2}-F_N
r^{-N/2+\sqrt{N-1}} ,\ \ \forall r\in (0,1]$, where $E_N>0$ and
$F_N$ only depend on $N$.
\end{enumerate}

Since $h$ is positive and increasing, then $\Phi''>0$. Hence
$(N-1)r^{N-3}u_r^2$ is increasing and so is $r^{2N-2}u_r^2$. This
implies that $-r^{N-1}u_r$ is increasing, which is is equivalent
to the positiveness of $g$.

On the other hand note that, since $0\leq h\leq 1$, we obtain
$\Phi'(r)\leq G_N r^{2\sqrt{N-1}-1}$ in $(0,1]$, for a constant
$G_N$. Hence, from the definition of $u_r$ we obtain ii).

To prove iii) observe that, from the positiveness of $h$, we
obtain $\Phi''(r)\geq r^{2\sqrt{N-1}}h'(r)$ in $(0,1]$. On the
other hand, from the definition of $u_r$ we have
$\Phi''(r)=(N-1)\left( (N-3)r^{N-4}u_r^2+2u_r u_{rr}
r^{N-3}\right)$. Therefore, by ii) and the previous inequality we
obtain iii).

Finally, it is easily seen that for every sequences
$\{r_n\}\subset(0,1]$, $\{y_n\}\subset \real^+$, with
$r_n\downarrow 0$, there exists $h\in C^2(0,1]$, increasing,
satisfying $0\leq h\leq 1$ and $h'(r_n)=y_n$. Take $y_n$ such that
$E_N y_n r_n^{-N/2+\sqrt{N-1}+2}-F_N r_n^{-N/2+\sqrt{N-1}}=M_n$.
Applying iii) we deduce $-u_{rr}(r_n)\geq M_n$ and the proof is
complete. \qed

\begin{corollary}\label{nohayy}

Let $N\geq 10$. There does not exist a function $\psi:(0,1]
\rightarrow \real^+$ with the following property: for every $u\in
H^1(B_1)$ semi-stable radially decreasing solution of a problem of
the type (\ref{mainequation}) with $g\geq 0$, there exist $C>0$
and $\varepsilon \in (0,1]$ such that $\vert u_{rr}(r)\vert \leq C
\psi (r)$ for $r\in (0,\varepsilon]$.

\end{corollary}

\noindent {\bf Proof.} Arguing as in Corollary \ref{nohay} and
using Proposition \ref{sucesiones2}, we conclude the proof of the
corollary. \qed

\begin{proposition}\label{sucesiones3}

Let $\{r_n\}\subset(0,1]$, $\{M_n\}\subset \real^+$ two sequences
with $r_n\downarrow 0$. Then, for $N\geq 10$, there exists $u\in
H^1(B_1)$, which is a semi-stable radially decreasing unbounded
solution of a problem of the type (\ref{mainequation}) with
$g,g'\geq 0$, satisfying

$$\vert u_{rrr}(r_n)\vert \geq M_n \ \ \ \ \forall n\in \mathbb{N}.$$

\end{proposition}

\begin{lemma}\label{hpequenna}

For any dimension $N\geq 10$, there exists $\varepsilon_N>0$ with
the following property: for every $h\in C^2(0,1]\cap C^1[0,1]$
satisfying $h(0)=0$, $0\leq h'\leq \varepsilon_N$ and $h''\leq 0$,
$u$ is a semi-stable radially decreasing unbounded $H^1(B_1)$
solution of a problem of the type (\ref{mainequation}) with
$g,g'\geq 0$, where $u_r$ is defined in terms of $h$ as in Theorem
\ref{familia}.

\end{lemma}

\noindent {\bf Proof.} Similarly as in the proof of Proposition
\ref{sucesiones2} (item i)), $h'\geq 0$ implies that $u$ is a
semi-stable radially decreasing unbounded $H^1(B_1)$ solution of a
problem of the type (\ref{mainequation}) with $g\geq 0$.

On the other hand, from the definition of $\Phi$ and $u_r$ it
follows easily that

$$\begin{array}{ll} u_r &\displaystyle{=-\sqrt{(N-1)^{-1}\,
r^{3-N}\Phi'}} \\ &\displaystyle{=
-r^{-N/2+\sqrt{N-1}+1}\sqrt{2(N-1)^{-1/2}\left(1+\int_0^r
h\right)+(N-1)^{-1}r\, h}}\\
\end{array}$$

Put this last expression in the form
$u_r=-r^{-N/2+\sqrt{N-1}+1}\varphi (r)$, where $\varphi(r)$ (and
of course $u_r$) depends on $h$. Now consider the set $X=\{ h\in
C^2(0,1]\cap C^1[0,1]:h(0)=0\, , 0\leq h'\, ,h''\leq 0 \}$ and the
norm $\Vert h\Vert_X=\Vert h'\Vert_{L^\infty (0,1)}$. Taking
$\Vert h\Vert_X \to 0$, we have

\begin{equation}\label{tyu}
\lim_{\Vert h\Vert_X \to 0} \varphi =\sqrt{2}(N-1)^{-1/4}, \ \
\lim_{\Vert h\Vert_X \to 0} \varphi'=0, \ \ \lim_{\Vert h\Vert_X
\to 0} \left( \varphi''-\frac{(N-1)^{-1}r\,
h''}{2\varphi}\right)=0,
\end{equation}

\noindent where all the limits are taken uniformly in $r\in
(0,1]$. On the other hand, it is easy to check that

$$\begin{array}{ll} r^2g'(u)&\displaystyle{=\frac{-r^2
u_{rrr}}{u_r}-\frac{(N-1)r\, u_{rr}}{u_r}+(N-1)}\\
&\displaystyle{=\frac{-r^2
\varphi''}{\varphi}-\frac{(2\sqrt{N-1}+1)\, r
\varphi'}{\varphi}+\frac{(N-2)^2}{4}}\\ \end{array}$$

Hence, from (\ref{tyu}), we can assert that, for $h\in X$ with
small $\Vert h\Vert_X$, $r^2 g'(u) >0$ in $(0,1]$, and the lemma
follows. \qed

\

\noindent {\bf Proof of Proposition \ref{sucesiones3}.} We follow
the notation used in the previous lemma. From (\ref{tyu}), we
deduce that

$$\lim_{\Vert h\Vert_X
\to 0} \left( r^{N/2-\sqrt{N-1}+1}u_{rrr}+\frac{(N-1)^{-1}\, r^3
h''}{2\varphi}\right)=\sigma ,$$

\noindent uniformly in $r\in (0,1]$, where $\sigma=-(
-N/2+\sqrt{N-1}+1)(-N/2+\sqrt{N-1})\sqrt{2}(N-1)^{-1/4}<0$. Then,
taking $\varepsilon'_N>0$ sufficient small (possibly less than
$\varepsilon_N$), we have that

$$r^{N/2-\sqrt{N-1}+1}u_{rrr}\geq -\left(\frac{(N-1)^{-1}\, r^3
h''}{2\sqrt{2}(N-1)^{-1/4}+1}\right)+\sigma -1\, , \ \forall r\in
(0,1],$$

\noindent for $\Vert h\Vert_X \leq \varepsilon'_N$.

Finally, it is easily seen that for every sequences
$\{r_n\}\subset(0,1]$, $\{y_n\}\subset \real^+$, with
$r_n\downarrow 0$, there exists $h\in X$, with $\Vert h\Vert_X
\leq \varepsilon'_N$, satisfying  $h''(r_n)=-y_n$. (Take, for
instance $h(r)=\int_0^r z(t)\, dt$, where $z\in C^1(0,1]\cap
C[0,1]$ is decreasing, $0\leq z(t)\leq \varepsilon'_N$ and
satisfies $z'(r_n)=-y_n$.) Take $y_n$ such that
$r_n^{N/2-\sqrt{N-1}+1}M_n=\left(\frac{(N-1)^{-1}\, r_n^3
y_n}{2\sqrt{2}(N-1)^{-1/4}+1}\right)+\sigma -1$. Applying the
above inequality, we obtain $u_{rrr}(r_n)\geq M_n$ and the proof
is complete. \qed

\begin{corollary}\label{nohayyy}

Let $N\geq 10$. There does not exist a function $\psi:(0,1]
\rightarrow \real^+$ with the following property: for every $u\in
H^1(B_1)$ semi-stable radially decreasing solution of a problem of
the type (\ref{mainequation}) with $g,g'\geq 0$, there exist $C>0$
and $\varepsilon \in (0,1]$ such that $\vert u_{rrr}(r)\vert \leq
C \psi (r)$ for $r\in (0,\varepsilon]$.

\end{corollary}

\noindent {\bf Proof.} Applying Proposition \ref{sucesiones3},
this follows by the same method as in Corollaries \ref{nohay} and
\ref{nohayy}. \qed

\

{\bf Acknowledgments.} The author would like to thank Xavier
Cabr\'e for very stimulating discussions.

\end{document}